\documentclass[final,leqno,letterpaper]{etna}

\setbibdata{1}{xx}{46}{2017} %

\usepackage{amsmath}
\usepackage{amsfonts}
\usepackage{graphicx}
\usepackage{epstopdf}
\usepackage{algorithmic}
\usepackage{bm}
\usepackage[capitalise, noabbrev]{cleveref} %
\usepackage{siunitx}
\usepackage{amsopn}
\usepackage{longtable}
\usepackage{multirow}
\usepackage{booktabs}

\usepackage{tikz}
\definecolor{tab-blue}{HTML}{1f77b4}
\definecolor{tab-orange}{HTML}{ff7f0e}
\definecolor{tab-green}{HTML}{2ca02c}
\definecolor{tab-red}{HTML}{d62728}
\definecolor{tab-purple}{HTML}{9467bd}
\definecolor{tab-brown}{HTML}{8c564b}
\definecolor{tab-pink}{HTML}{e377c2}
\definecolor{tab-gray}{HTML}{7f7f7f}
\definecolor{tab-olive}{HTML}{bcbd22}
\definecolor{tab-cyan}{HTML}{17becf}

\tikzset{every picture/.style={>=latex,
                               font={\fontsize{9}{11pt}\selectfont},
                               line width=2pt,
                               line cap=round,
                              },
         mainaxis/.style={line width=0.6pt},
         tickmajor/.style={line width=0.6pt},
         tickminor/.style={line width=0.4pt},
         graphnode/.style={draw=tab-blue,
                           line width=2pt,
                           circle,
                           align=center},
         graphedge/.style={draw=tab-gray,
                           line width=1pt,
                          }
        }

\newcommand{\thelongtitle}{Unstructured to structured: geometric multigrid on complex geometries via domain remapping}
\newcommand{\theshorttitle}{Geometric multigrid via domain remapping}

\hypersetup{%
    pdftitle={\theshorttitle},
    pdfauthor={Nicolas Nytko, Scott MacLachlan, J. David Moulton, Luke Olson, Andrew Reisner, Matt West},
    pdfkeywords={geometric multigrid, boxmg, diffeomorphism, learning, neural ode}
    }

\title{\thelongtitle\thanks{%
Received... Accepted... Published online on... Recommended by....
}}

\author{Nicolas Nytko\footnotemark[2]
  \and Scott MacLachlan\footnotemark[3]
  \and J. David Moulton\footnotemark[4]
  \and Luke N. Olson\footnotemark[2]
  \and Andrew Reisner\footnotemark[5]
  \and Matthew West\footnotemark[6]}

\shorttitle{\theshorttitle}
\shortauthor{Nytko, MacLachlan, Moulton, Olson, Reisner, West}

\newcommand{\dx}{\ dx}

\newcommand{\Div}{\nabla \cdot}
\newcommand{\Grad}{\nabla}
\newcommand{\mat}[1]{\bm{#1}}
\newcommand{\tens}[1]{\bm{\mathcal{#1}}}
\renewcommand{\vec}[1]{\bm{#1}}
\newcommand{\x}{\vec{x}}
\newcommand{\Vp}{\mathcal{V}_p}
\newcommand{\Vc}{\mathcal{V}_c}

\newcommand{\D}{\mathcal{D}}
\newcommand{\Dc}{\D_c}

\DeclareMathOperator{\cond}{cond}
\DeclareMathOperator{\tr}{tr}
\DeclareMathOperator{\nnz}{nnz}
\DeclareMathOperator{\atantwo}{atan2}
\DeclareMathOperator*{\argmin}{argmin}
\begin{document}

\maketitle

\renewcommand{\thefootnote}{\fnsymbol{footnote}}

\footnotetext[2]{Siebel School of Computing and Data Science, University of Illinois Urbana-Champaign}
\footnotetext[3]{Mathematics and Statistics, Memorial University of Newfoundland}
\footnotetext[4]{Energy and Natural Resources Security (EES-16), Los Alamos National Laboratory}
\footnotetext[5]{Computer, Computational, and Statistical Sciences (CCS-7), Los Alamos National Laboratory}
\footnotetext[6]{Mechanical Sciences and Engineering, University of Illinois Urbana-Champaign}

\begin{abstract}
For domains that are easily represented by structured meshes, robust geometric multigrid solvers
can quickly provide the numerical solution to many discretized elliptic PDEs. However, for complicated domains with unstructured meshes,
constructing suitable hierarchies of meshes becomes challenging.  We propose a framework for mapping
computations from such complex domains to regular computational domains via diffeomorphisms,
enabling the use of robust geometric-style multigrid.  This mapping facilitates regular memory accesses during
solves, improving efficiency and scalability, especially on massively parallel processors such as
GPUs.  Moreover, we show that the diffeomorphic mapping itself may be approximately learned using an
invertible neural network, facilitating automated application to geometries where no analytic mapping
is readily available.
\end{abstract}



\section{Introduction}

Geometric multigrid methods are efficient iterative solvers for elliptic partial differential
equations (PDEs), offering an optimal time complexity and rapid convergence for problems where
geometric structure can be exploited.  However, extending multigrid to arbitrary meshes --- finite
partitions of domains, which are open subsets of $\mathbb{R}^{d}$ --- remains challenging due to the
difficulty of constructing a suitable hierarchy of coarse meshes that maintains multigrid
efficiency.  In practice, algebraic multigrid (AMG) methods~\cite{Ruge_1987,KStuben_2001a} are often used on
complex meshes as they do not require geometric information about the problem, operating directly on
the linear system itself.  Despite their versatility, AMG approaches tend to suffer from higher
setup costs, as more expensive graph algorithms are required to coarsen the problem.  Moreover, AMG can
suffer from poor data locality and irregular data accesses due to sparse-matrix storage patterns
that hinder performance on massively parallel architectures such as
GPUs~\cite{Bell2012}.

To address these limitations, we propose a novel framework to expand the applicability of
geometric multigrid methods to a larger class of problems.  To do so, we construct diffeomorphic mappings between a given mesh and a structured (regular) mesh on which robust geometric multigrid methods can be used.  A full
multigrid hierarchy is then constructed by transferring solution and/or residual vectors from the given (unstructured) mesh to
the structured one and employing a black-box multigrid
solver~\cite{DendyBoxMG1982,DendyBoxMG2010,ReisnerBoxMG2018} to quickly and efficiently solve the
structured problem.  This auxiliary mesh correction or solution is then transferred back to the unstructured mesh, where
only cheap, local relaxation is performed.  This results in a solver that leverages the regularity of
the structured problem to yield an efficient solver, while also allowing for essential features of
the physical geometry to be preserved.

Additionally, we show that the diffeomorphic mapping between domains can be
learned.  Taking inspiration from existing work on large deformation diffeomorphic metric
mappings~\cite{Hernandez2024} (LDDMMs), used heavily in computational anatomy, we demonstrate that
such a mapping can be approximately learned for geometries lacking an analytic conformal mapping.
We show numerical results for both domains with existing analytic mappings and ones where the
mapping itself is learned; iteration counts to convergence and time-to-solution of the resulting geometric solver are
compared to algebraic methods, showing that our framework is relatively robust and comparable with
other multilevel solvers.

A variety of related works have explored the usage of structured discretizations for numerical PDE
problems posed on complex geometries, both using the structure to outright solve the PDE and also to
form a nested multigrid hierarchy for efficient computation.  Cut finite element methods
(CutFEM)~\cite{BurmanCutFEM2015} overlay a structured background mesh over a geometric description
of a domain, weakly enforcing boundary and interface conditions through Nitsche's
method~\cite{JuntunenNitsche2009} and additional penalty terms to enforce regularity.  Overlapping
grid methods~\cite{HenshawOverlapping2003,ChesshireOverlapping1990,HenshawOverlapping2005}
approximate a complex ``Chimera'' domain by combining multiple overset structured grids with simple
(shear) transforms applied.  There also exist works studying projection-based interpolation of
functions between nonmatching overlapping finite element
spaces~\cite{PontInterpolation2016,XiangminNonmatching2004}, optionally imposing restrictions such
as weak conservation of mass or preserving boundary norms; this is, in effect, similar to how we
compute the projection operator for our multilevel solver.  In multigrid contexts, the auxiliary
space method~\cite{XuAuxiliary1996,ChenAuxiliary2015,BrambleNonnested1991} constructs a structured
auxiliary problem and stable transfer operators to precondition unstructured (fine level)
discretizations; in constrast, our framework builds this auxiliary space through smooth geometric
mappings.  Recently in the scientific machine learning domain, diffeomorphic transformations have
been similarly used to solve PDEs on regular domains by use of Fourier Neural Operators
(FNOs)~\cite{LiFNO2023,LiFNO2021}; similarly to this work, both analytic and learned mappings are
considered.

Throughout this work, we will explicitly use the term \emph{domain} to refer to an open subset of
$\mathbb{R}^d$ and \emph{mesh} to refer to its respective finite partition composed of simplices or
tensor-product elements.  This distinction helps distinguish between the continuous geometric setting
and the discrete numerical approximation.

The remainder of this paper is outlined as follows.  We introduce our method in~\cref{sec:method},
including details on the multigrid method in~\cref{ssec:multigrid} and the grid-transfer operators
between the structured and unstructured grids in~\cref{subsec:interpolation,subsec:cgc_scaling}.
Details of the construction and learning of mappings between structured and unstructured grids is
presented in~\cref{sec:mappings}.  Numerical results are given in~\cref{sec:numerics}, with
conclusions given in~\cref{sec:conclusions}.

\section{Method}\label{sec:method}
We consider the elliptic diffusion problem with Dirichlet boundary
conditions,
\begin{align}
  -\Div \D \Grad u &= f \quad \text{ in } \Omega_p, \label{eqn:pde} \\
  u &= g \quad \text{ on } \partial \Omega_p,
\end{align}
defined over an open, connected space (the ``physical domain'') $\Omega_p \subset \mathbb{R}^d$ with continuous boundary
$\partial \Omega_p$ and symmetric, positive-definite diffusion tensor $\D : \mathbb{R}^d \to
\mathbb{R}^{d \times d}$.  We assume there exists a triangulation $\mathcal{T}_p = \big\{ \tau_p^i
\big\}$ such that $\Omega_p = \cup_i \tau_p^i$, although we will demonstrate below that our approach naturally generalizes to the case where $\Omega_p \approx \cup_i \tau_p^i$, such as when the physical domain has curved boundaries.  In this work, we consider the case of $d=2$, though the method directly
extends to $d=3$ as well.  Similarly, we focus on the diffusion equation given in~\cref{eqn:pde}, although the method could also be applied to reaction-diffusion or convection-diffusion equations.

\begin{figure}[!htp]
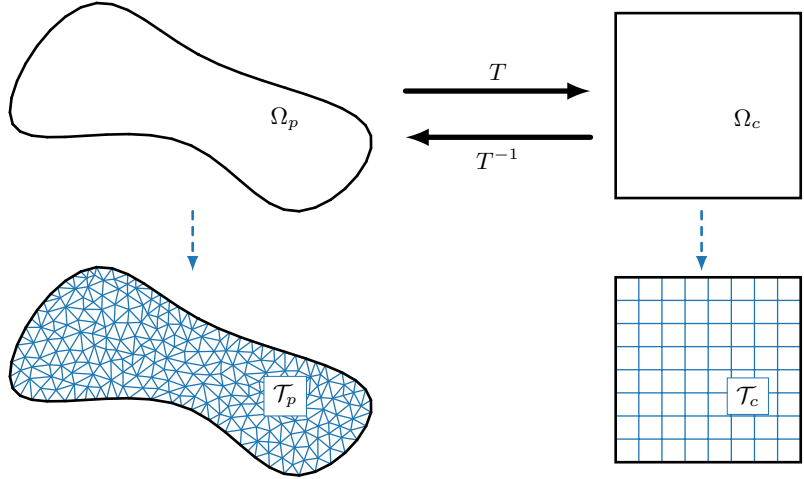

  \centering


  \caption{Schematic of the \emph{physical} domain (upper left), the \emph{computational} domain (upper right), and their respective meshes, $\mathcal{T}_{p}$ and $\mathcal{T}_c$.  The mapping $T:\,\Omega_p \rightarrow \Omega_c$ is assumed invertible.}\label{fig:domain}
\end{figure}

In addition to the physical domain, $\Omega_p$, we define an auxiliary \emph{computational domain},
$\Omega_c = (-1,1)^d$ and, as discussed below, a corresponding triangulation, $\mathcal{T}_c =
\big\{ \tau_c^i \big\}$ with $\Omega_c = \cup_i \tau_c^i$.  As result, both $\Omega_p$ and
$\Omega_c$ are smooth manifolds, allowing us to additionally define an invertible map $T :
\overline{\Omega}_p \mapsto \overline{\Omega}_c$ between the two domains, including their
boundaries.  Moreover, we require $T$ to be smooth (diffeomorphic) over the interior to ensure the
Jacobian $\mat{J}_T(\x)$ and its inverse exist and are continuous.  We use $|\cdot|$ to denote the absolute value of the determinant when applied to a matrix; for example, $|\mat{J}_T| = |\det(\mat{J}_T)|$.  In what follows, we will assume
$\mathcal{T}_c$ corresponds to a logically structured grid of regular elements, such as a
tensor-product mesh of quadrilateral elements.  A schematic of the domains and meshes is shown
in~\cref{fig:domain}.  We first consider example domains where $T$ is known analytically. In
\cref{subsec:learning}, however, we will show that $T$ can also be learned for more complex domains.

Our method assumes knowledge of both the continuum PDE, given in~\cref{eqn:pde}, and of the mapping,
$T$.  As shown below, this allows us to define an equivalent weak form of the PDE over the
computational domain.  We make use of this to define a structured multigrid hierarchy (using robust
geometric multigrid components) on which most of the computational work in a multigrid V-cycle is
performed.  It is well known that structured-grid methods are generally the most efficient and
straightforward of multigrid methods, due to their highly regular memory access patterns that do not
need the indirection typical of generic sparse matrix formats, such as compressed sparse row (CSR)
format.  This clearly has high potential to yield speedups on compute platforms, such as GPUs, that
have performance bottlenecks resulting from memory transfers, rather than raw compute power.

We begin with the regular weak form for the diffusion equation on the physical domain, $\Omega_p$,
where we write the solution to~\cref{eqn:pde} as $u(\vec{x}) = \hat{u}_p(\vec{x}) +
\hat{g}(\vec{x})$, where $\hat{g} \in H^1(\Omega_p)$ satisfies the boundary condition
$\hat{g}(\vec{x}) = g$ on $\partial \Omega_p$, and seek the function $\hat{u}_p \in H^1_0(\Omega_p)$
that satisfies
\begin{equation}
  \int_{\Omega_p} \D \Grad \hat{u}_p \cdot \Grad v \dx = \int_{\Omega_p} f v \dx - \int_{\Omega_p} \D \Grad \hat{g} \cdot \Grad v \dx \quad\quad
  \forall v \in H^1_0(\Omega_p), \label{eqn:weak_physical}
\end{equation}
where $H^1(\Omega_p)$ is the standard Sobolev space of square-integrable functions
over $\Omega_p$ with square-integrable first derivatives, and $H^1_0(\Omega_p)$ is the restriction of
$H^1(\Omega_p)$ to functions that vanish on the boundary.
The weak form can be mapped to the
computational domain by the change of coordinates defined by $T$, in a similar
fashion to standard finite-element assembly on a general, triangular mesh.
We transform the weak form in~\cref{eqn:weak_physical} to one over
$\Omega_c$ using $\mat{J}_T$, asking to find $\hat{u}_c \in H^1_0(\Omega_c)$ such that
\begin{align}
  \int_{\Omega_c} ((\D \circ T^{-1}) &\mat{J}_T^T \Grad \hat{u}_c) \cdot (\mat{J}_T^T \Grad v)  \ |\mat{J}_T^{-1}| \dx =
\int_{\Omega_c} (f \circ T^{-1}) v \ |\mat{J}_T^{-1}| \dx \notag \\ & - \int_{\Omega_c} \left((\D\circ T^{-1}) \mat{J}_T^T \Grad (\hat{g}\circ T^{-1})\right) \cdot (\mat{J}_T^T \Grad v) |\mat{J}_T^{-1}|\dx \quad\quad \forall v \in H^1_0(\Omega_c). \label{eqn:weak_computational}
\end{align}
We note that we write $\hat{u}_p$ for the solution in the physical domain and $\hat{u}_c$ for the
solution in the computational domain, which are related by composition with $T$ or $T^{-1}$, while
we explicitly compose the right-hand side data, $f$ and $\hat{g}$ (defined on the physical domain),
with $T^{-1}: \Omega_c \rightarrow \Omega_p$ when writing the weak form in the computational domain.
In practice, we will not need to evaluate the right-hand side in~\cref{eqn:weak_computational}, so
we need not consider the details of these mappings.

\subsection{Discretization}\label{subsec:discretization}
We discretize the weak forms over the physical and computational domains,
\cref{eqn:weak_physical,eqn:weak_computational}, by restricting the test and trial function spaces
to discrete function spaces $\Vp \subset H^1_0(\Omega_p)$ and $\Vc \subset H^1_0(\Omega_c)$.  As is natural with finite-element methods, these function spaces are associated with the meshes, $\mathcal{T}_p$ and $\mathcal{T}_c$, although we note that both the meshes and the discretization spaces can be chosen completely independently, as depicted in~\cref{fig:domain}.  For both $\Vp$ and $\Vc$, we assume there exist known bases such that any function, $f_p \in \Vp$ can be written as a weighted sum of finitely many basis functions
$f_p = \sum_i^{N_p} f^i_p \phi^i_p$, and with a similar expression for any $f_c \in \Vc$.
We then write $\vec{f}_p \in \mathbb{R}^{N_p}$ as the vector whose entries are the
coefficients $\{f_p^i\}$.

While the framework we present is general enough for any discrete bases to be used for the physical and
computational spaces, in this work we consider the case where $\mathcal{V}_p = P_1(\mathcal{T}_p)$, piecewise linear functions on a triangular mesh of the physical domain, and $\mathcal{V}_c = Q_1(\mathcal{T}_c)$, piecewise bilinear functions on a tensor-product quadrilateral mesh of the computational domain.  The
integrals in the bilinear forms from~\cref{eqn:weak_physical,eqn:weak_computational} are numerically approximated and assembled into matrices, $\mat{A}_p$ and $\mat{A}_c$, for the physical and computational operators,
respectively.  This integration and assembly step is done as standard for finite elements; in
practice, the order required of the quadrature rule depends on the smoothness of $\mathcal{D}$ and the regularity of the mapping $T$,
with rapidly changing PDE coefficients and/or ill-conditioned maps (with rapid oscillations/changes) requiring higher-order quadrature schemes.  Thus, when constructing maps to the computational domain, we aim to compute maps that are as smooth as possible, to minimize the induced computational cost of the mapping.

An important consequence of the mapping to the computational domain is that we must solve a ``harder'' problem there than on the physical domain, due to the mapping.  One way to quantify the added difficulty is to consider the effective diffusion coefficient, $\Dc$, on the computational domain, given by
  \begin{equation}
    \Dc = \mat{J}_T \big(\D \circ T^{-1}\big) \mat{J}_T^T |\mat{J}_T^{-1}| .
  \end{equation}
First of all, we note that there is a natural (but not sharp) bound on the condition number of $\Dc$, $\cond(\Dc)$, as
\begin{equation}
  \cond(\Dc(\hat{\x})) \leq \cond(\mat{J}_T)^2 \cond(\D(\x)),
\end{equation}
when $\hat{\x} = T^{-1}(\x)$ for some $\x \in \Omega_p$.  Thus, $\Dc$ reflects any anisotropy present in $\D$, while ill-conditioning present in $T$ can be amplified quadratically in the resulting diffusivity coefficient.  This can be either helpful or harmful, depending on whether these align with those of $\D$ and, if they do, whether they cancel out ill-conditioning in $\D$ or compound it.  Secondly, we note that rewriting the diffusivity coefficient in this form allows easy use of existing finite-element
discretizations on structured grids, without needing special treatment of the weak form in
\cref{eqn:weak_computational}.

\subsection{Multigrid Solver}\label{ssec:multigrid}
Our solution methodology makes use of the following ingredients:
\begin{enumerate}
\item The linear system for the PDE on the physical domain, $\mat{A}_p \vec{u}_p = \vec{f}_p$,
 either defined through an explicit matrix, $\mat{A}_p$, or by way of a method to compute matrix-vector products;
\item a stationary iteration, $S$, to (cheaply) relax residuals of $\mat{A}_p\vec{u}_p = \vec{f}_p$;
\item Matrix $\mat{A}_c$ corresponding to the bilinear form on the left-hand side of~\cref{eqn:weak_computational}; and
\item an interpolation operator, $\mat{P}$, mapping from $\mathcal{V}_c$ to $\mathcal{V}_p$ as described below in~\cref{subsec:interpolation}.
\end{enumerate}

\begin{figure}[h]
  \centering
  \includegraphics{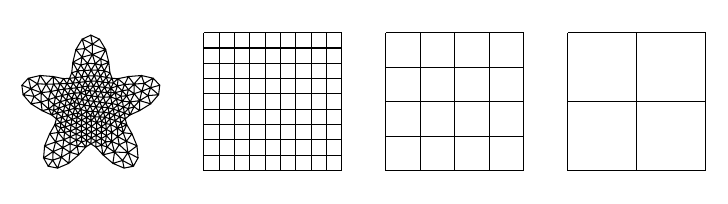}
  \caption{An example multigrid hierarchy.  The finest level (left-most mesh) corresponds to the mesh for the physical domain, followed by the hierarchy of structured grids (right three
meshes) on the computational domain.}\label{fig:mg_hierarchy}
\end{figure}

We solve the physical system, $\mat{A}_p\vec{u}_p = \vec{f}_p$, using a two-grid multigrid method,
with the computational grid playing the role of the coarse-grid operator, and $\mat{P}$ and
$\mat{P}^T$ as the grid-transfer operators.  We solve the coarse-grid system, with matrix
$\mat{A}_c$, by applying the black-box multigrid (BoxMG)
algorithm~\cite{DendyBoxMG1982,DendyBoxMG2010} as implemented in the Cedar software
package~\cite{ReisnerBoxMG2018}, though any structured-grid multigrid solver may be used.  The
method in Cedar makes use of structured geometric coarsening (by a factor of two in each direction)
and operator-induced interpolation, giving robust solve capabilities when anisotropy is present in
the problem; this is combined with Galerkin coarsening to generate the coarse grids.  Thus, the full
hierarchy in our approach is formed by prepending the original physical system to the computational
hierarchy, as shown in~\cref{fig:mg_hierarchy}.  We present details of the solver algorithm in the
following subsections.

\subsubsection{Setup phase}  The setup phase of the algorithm is primarily that of BoxMG on the matrix $\mat{A}_c$, but also includes setting up the transfers between the finest computational grid and the physical grid.

As inputs, the setup phase requires knowledge of the physical and computational meshes, $\mathcal{T}_p$ and $\mathcal{T}_c$, as well as the mapping, $T$.  Either $\mat{A}_c$ must be provided, or the diffusion coefficient, $\D$, must be known.  With this, the setup phase of the algorithm requires:
\begin{enumerate}
\item Compute interpolation $\mat{P}$ by mapping vertices of $\mathcal{T}_p$ via $T$
(see \cref{subsec:interpolation}).
\item If needed, compute $\mat{A}_c$ from basis for $\mathcal{V}_c$ and~\cref{eqn:weak_computational}.
\item Setup BoxMG solver on $\mat{A}_c$.
\end{enumerate}

\subsubsection{Solution phase}  The solution phase of the algorithm consists of stationary iterations wrapped around a V-cycle multigrid structure.  This requires the physical grid approximate solution, $\vec{u}_p$, right-hand side, $\vec{f}_p$, and matrix, $\mat{A}_p$, as well as the relaxation operator, $S$, and the computational grid BoxMG solver.  A single V-cycle then takes essentially the standard form:
\begin{enumerate}
\item Pre-relaxation: $\vec{u}_p \gets S(\vec{f}_p - \mat{A}_p \vec{u}_p)$.
\item Restrict residual: $\vec{r}_c \gets \mat{P}^T (\vec{f}_p - \mat{A}_p \vec{u}_p)$.
\item BoxMG V-cycle on $\mat{A}_c\vec{u}_c = \vec{r}_c$.
\item Prolong correction: $\vec{u}_p \gets \vec{u}_p + \mat{P} \vec{u}_c$.
\item Post-relaxation: $\vec{u}_p \gets S(\vec{f}_p - \mat{A}_p \vec{u}_p)$.
\end{enumerate}

\subsection{Interpolation}\label{subsec:interpolation}

Clearly the success or failure of the proposed method depends strongly on the quality of the grid-transfer operator, $\mat{P}$, from $\Omega_c$ to $\Omega_p$.  We define this based on the
minimization problem: given a function, $u_c \in \Vc$, we want to find the \emph{nearest} (in an $L^2$ sense)
function
$u_p \in \Vp$ under the map $T$:
\begin{equation}
u_p  = \argmin_{u_p \in \Vp} \frac{1}{2} \| u_p - u_c \circ
T\|^2_{L^2(\Omega_p)}. \label{eqn:argmin_functional}
\end{equation}
This is minimized by solving the associated Euler-Lagrange equation~\cite{GelfandCalculus2000},
leading to the equality
\begin{equation}  \int_{\Omega_p} u_p v \dx = \int_{\Omega_p} (u_c \circ T) v \dx \quad \forall v
\in \Vp, \label{eqn:functional_min}
\end{equation}
or, equivalently, when $u_c$ and $u_p$ are expanded out in terms of their respective basis functions,
\begin{equation}
\sum_{j=1}^{N_p} u_p^j
\int_{\Omega_p} \phi_p^j \phi^i_p \dx =   \sum_{j=1}^{N_c} u_c^j \int_{\Omega_p} (\phi^j_c \circ T) \phi_p^i \dx \quad \forall \phi_p^i \in \Vp.
\end{equation}

Since the equation is linear in the coefficients $u_p^j$ and $u_c^j$, we
can express the above equality as
\begin{equation}
  \mat{M}_p \vec{u}_p = \mat{C} \vec{u}_c \implies \vec{u}_p = \mat{M}_p^{-1} \mat{C} \vec{u}_c.
\end{equation}
where the bilinear \emph{physical} and \emph{coupled} mass forms, respectively, are defined by
\begin{align}
  \big[ \mat{M}_p \big]_{ij} &= \int_{\Omega_p} \phi^j_p \phi^i_p \dx, \label{eqn:physical_mass} \\
  \big[ \mat{C} \big]_{ij} &= \int_{\Omega_p} (\phi^j_c \circ T) \phi^i_p  \dx. \label{eqn:coupled_mass}
\end{align}
As $\vec{u}_c$ is an arbitrary function in $\Vc$, the general prolongation operator
becomes $\mat{P} = \mat{M}_p^{-1} \mat{C}$.  In practice, when we apply $\mat{P}$ to a vector, we only approximately invert the mass matrix, using just 2 iterations of unpreconditioned CG.

When seeking a restriction from the physical to the computational space, we aim to find the
functional $u_c^\ast \in \Vc^\ast$ that is nearest to a given $u_p^\ast \in \Vp^\ast$.  An important observation is
that we are seeking elements in the dual of the function space; in a multigrid setting, residuals are
defined by a functional corresponding to evaluations of the bilinear form and the right-hand side.

By the Riesz representation theorem, we have unique $\hat{u}_c, \hat{u}_p$ such that
\begin{align}
  u_c^\ast(v) &= \int_{\Omega_c} \hat{u}_c v \dx \quad \forall v \in \Vc, \\
  u_p^\ast(v \circ T) &= \int_{\Omega_p} \hat{u}_p (v \circ T) \dx = \int_{\Omega_c} (\hat{u}_p
\circ T^{-1}) v \ |\mat{J}_T^{-1}|\dx \quad \forall v \in \Vc,
\end{align}
thus allowing us to set up a similar minimization problem to \cref{eqn:argmin_functional}:
\begin{equation}
  \begin{split}
  u_c^\ast = &\argmin_{u_c^\ast \in \Vc^\ast} \frac{1}{2} \big\|u_c^\ast - u_p^\ast \circ
T\big\|^2_{(L^2(\Omega_c))^\ast} \\
  &\quad= \argmin_{u_c^\ast \in \Vc^\ast} \frac{1}{2} \big\|\hat{u}_c - (\hat{u}_p \circ T^{-1})
|\mat{J}_T|^{-1}\big\|^2_{L^2(\Omega_c)}.
  \end{split}
\end{equation}
Again solving the associated Euler-Lagrange equation and expanding in terms of basis coefficients
gives
\begin{align*}
  \sum_{j=1}^{N_c} \hat{u}_c^j \int_{\Omega_c} \phi_c^j \phi_c^i \dx
  &= \sum_{j=1}^{N_p} \hat{u}_p^j \int_{\Omega_c} (\phi_p^j \circ T^{-1}) \phi_c^i |\mat{J}_T|^{-1}\dx
\quad \forall \phi_c^i \in \Vc \\
  &= \sum_{j=1}^{N_p} \hat{u}_p^j \int_{\Omega_p} \phi_p^j (\phi_c^i \circ T)\dx \quad \forall
\phi_c^i \in \Vc,
\end{align*}
which leads to a similar linear equality as above with the computational mass matrix $\mat{M}_c$,
\begin{equation}
  \mat{M}_c \hat{\vec{u}}_c = \mat{C}^T \hat{\vec{u}}_p. \label{eqn:restriction_fn}
\end{equation}
Interpreting the physical and computational mass matrices as Riesz maps on $\Vp$ and $\Vc$,
respectively, we have
\begin{equation}
  \mat{M}_p \hat{\vec{u}}_p = \vec{u}_p^\ast \quad
  \mat{M}_c \hat{\vec{u}}_c = \vec{u}_c^\ast,
\end{equation}
allowing us to rewrite \cref{eqn:restriction_fn} as
\begin{align}
   \vec{u}_c^\ast = \mat{C}^T \mat{M}_p^{-1} \vec{u}_p^\ast,
\end{align}
implying that restriction is equal to the transpose of interpolation, $\mat{R}=\mat{P}^T$.

\subsubsection{Coupled Mass Evaluation}
Computation of the coupled mass matrix, $\mat{C}$, requires evaluating basis functions on both the
physical and computational meshes and becomes slightly more involved when we do not have degrees of
freedom that align between the two meshes.  However, applying a change of variables to
\cref{eqn:coupled_mass} via $T$, we observe that there are, indeed, two equivalent formulations of the
weak form defining $\mat{C}$,
\begin{equation}
  \big[\mat{C}\big]_{ij} = \int_{\Omega_p} \phi^i_p (\phi^j_c \circ T) \dx =
  \int_{\Omega_c} (\phi^i_p \circ T^{-1}) \phi^j_c \ |\mat{J}_T^{-1}| \dx, \label{eqn:coupled_mass_comp}
\end{equation}
giving us some choice in how to evaluate the integral.

We elect to compute the integral in computational space, as our above assumption of using
uniformly spaced grids on the computational domain greatly reduces the required computational effort, as it can be written in terms of evaluation of
tensor-product elements on an axis-aligned grid.  The algorithm proceeds as follows: we iterate
through each element $\tau_p^i \in \mathcal{T}_p$ on the physical mesh and compute $\tilde{\tau}_p^i
\subset \Omega_c$ by mapping the \emph{vertices} of $\tau_p^i$ onto $\Omega_c$ through $T$.  From
this, we determine the subset $W_i = \big\{ \tau^j_c : \tau^j_c \cap \tilde{\tau}_p^i \neq \emptyset
\big\}$, of computational elements that overlap with $\tilde{\tau}_p^i$; in practice, this is done in
a process that is similar to rasterization in computer graphics.

From this overlapping set, we loop through pairs of overlapping elements on both domains,
$(\tau_p^i, \tau_c^j)$, and compute a partial contribution of the integral in
\cref{eqn:coupled_mass_comp}.  The quadrature points on $\tau_c^j$ are mapped to a reference element
on the physical mesh, $\hat{\tau}_p$.  Quadrature points that lie outside of the reference element
(which correspond to regions where the elements do not intersect when there is only partial overlap)
are simply evaluated as $0$, whereas points on the interior are used to evaluate the weak form as
usual.

By determining overlapping sets of elements, we evaluate $\mat{C}$ in a way that is stable when
one mesh is \emph{highly refined} compared to the other.  If we were to simply loop over elements on
one mesh and map their quadrature points to the other domain,
quadrature aliasing may occur, where certain elements do not receive any mass contribution even though
the elements themselves geometrically intersect.

In our implementation, we evaluate and store $\mat{C}$ as a sparse matrix.  However, we note that
an implementation variant (and perhaps, more friendly to GPUs) is to precompute the overlapping sets $W_i$
and compute the action of $\mat{C}$ on a vector via a matrix-free algorithm.

\subsection{Coarse-Grid Scaling}\label{subsec:cgc_scaling}

Preliminary numerical experiments showed that, when using low-order numerical quadrature to evaluate
$\mat{C}$, a significant ``mis-scaling'' could arise between the discretized problem on the
computational grid and the original one on the physical grid.  This leads to slow convergence, where
the \emph{direction} of the coarse-grid correction is a good one, but the scaling of the correction
is poorly chosen relative to the current approximation, $\vec{u}_p$.  This was observed to be
particularly prevalent when the map $T$ is ill-conditioned and has an oscillatory Jacobian.
Consequently, we propose a heuristic coarse-grid scaling that helps to correct the energy-mismatch
and, in many cases, accelerates convergence of the solver.

This scaling takes the simple form of under/over-correction based on the coarse-grid correction computed from $\mat{A}_c\vec{u}_c = \vec{r}_c$, taking $\vec{u}_p \gets \vec{u}_p + \gamma \mat{P}\vec{u}_c$, for scaling factor, $\gamma$.  Ideally, we would define $\gamma$ as
\begin{equation}
 \gamma = \argmin_{\alpha} \sum_{i=1}^k (\vec{x}_i^T \mat{P}^T \mat{A}_p \mat{P} \vec{x}_i - \alpha \vec{x}_i^T \mat{A}_c \vec{x}_i)^2, \label{eqn:cgc_scaling_1}
\end{equation}
over a set of test vectors $\vec{x}_i \in \big\{\vec{x}_1, \ldots, \vec{x}_k\big\}$, where these vectors represent pre-images of errors that are in the range of interpolation, $\mat{P}$ and are slow to converge under relaxation, $S$, since we solve the problem $\mat{A}_c \vec{u}_c = \mat{P}^T \mat{A}_p \vec{e}_p$ on the computational grid.  Thus, while
these test vectors could be sampled either adaptively (as in~\cite{MBrezina_etal_2006a}) or randomly, we choose them by considering the solution to the generalized eigenvalue problem
\begin{equation}
  \mat{P}^T \mat{A}_p \mat{P} \vec{x}_i = \lambda_i \mat{A}_c\vec{x}_i \quad \vec{x}_i^T \mat{A}_c \vec{x}_i = 1, \label{eqn:gen_eigvalue}
\end{equation}
to obtain an $\mat{A}_c$-orthogonal basis of vectors.  This greatly simplifies the optimization
problem in \cref{eqn:cgc_scaling_1}, allowing us to instead write
\begin{equation}
\gamma =   \argmin_{\alpha} \sum_{i=1}^k (\lambda_i \vec{x}_i^T\mat{A}_c\vec{x}_i - \alpha \vec{x}_i^T \mat{A}_c \vec{x}_i)^2,
\end{equation}
which, by assumption that $\vec{x}_i^T \mat{A}_c \vec{x}_i = 1$, gives
\begin{equation}
\gamma =  \argmin_{\alpha} \sum_{i=1}^k (\lambda_i - \alpha)^2 = \frac{1}{k} \sum_{i=1}^k \lambda_i,
\end{equation}
or that $\gamma$ is the arithmetic mean of the generalized eigenvalues.

In practice, computing all pairs $(\lambda_i, \vec{x}_i)$ to \cref{eqn:gen_eigvalue} is prohibitively
expensive as the size of the computational grid grows larger.  Experiments have shown that it is sufficient to sample a handful ($k=12$) of
eigenvalues through, for example, a generalized Lanczos iteration~\cite{SaadEigs2011}.

\section{Mappings}\label{sec:mappings}
The solver outlined in~\cref{ssec:multigrid} requires the mapping, $T$, from the physical to computational domain as input.  For many domains, such mappings to a square are known a priori.  However, for more
complicated geometries such a mapping is likely not available.  We first provide background theory on
mappings between domains, motivating the learning of the mapping discussed in~\cref{subsec:learning}.

\begin{definition}[Harmonic map]
  For some open, connected sets $U, V \subset \mathbb{R}^d$, a homeomorphism $T : \overline{U} \to
\overline{V}$ is said to be \emph{harmonic} if its components each individually satisfy Laplace's
equation on the interior, i.e.,
  \begin{equation}
    -\Delta T(\x) = \mat{0} \quad \forall \x \in U.
  \end{equation}
  Additionally, if we have a boundary homeomorphism $\phi : \partial U \to \partial V$ and let
  $T$ satisfy the boundary condition
  \begin{equation}
    T|_{\partial U} = \phi,\label{eq:boundary_homeomorphism}
  \end{equation}
  then we refer to $T$ as the \emph{harmonic extension} of $\phi$.
  \label{def:harmonic_map}
\end{definition}

Harmonic maps arise from minimization of the Dirichlet energy of the map, $T$,
\begin{equation}
  E_d[T] = \int_{\Omega} \|\nabla T\|^2_2 \dx,
\end{equation}
constrained by the boundary condition in~\cref{eq:boundary_homeomorphism}.
As they are critical points of the functional, they can be considered generalizations of geodesics to
maps between manifolds, in that they encode the shortest paths between points, given a distance
metric.  In this work, we prefer to use harmonic maps, as they tend to have smooth structure and
well-behaved derivatives.  Specifically, it is possible to construct a harmonic map, $T$, that is
diffeomorphic on $\Omega_p \to \Omega_c$ in two dimensions.  Here, we consider Theorem 4 from
\cite{Alessandrini2001}.

\begin{theorem}
  Let $U, V \subset \mathbb{R}^2$ be bounded, simply connected, open convex sets, and $\phi : \partial
  U \to \partial V$ a homeomorphism continuously mapping between the boundaries of $U$ and $V$.  If
  $T$ is a harmonic extension of $\phi$ satisfying
  \begin{align}
    -\Delta T &= \mat{0} \quad \text{in } U, \\
    T &= \phi \quad \text{on } \partial U,
  \end{align}
  then $T$ is a homeomorphism of $\overline{U}$ onto $\overline{V}$.  Moreover, if
  the Jacobian determinant is strictly positive, $|J_T(\x)|>0 \enskip \forall \x \in U$, then
  $T$ is smooth (and thus, diffeomorphic) on the interior.
  \label{thm:diffeo}
\end{theorem}

\begin{remark}
  From \cite{Alessandrini2001}, the boundary map $\phi$ being unimodal in each coordinate is a
sufficient condition to ensure the Jacobian determinant is strictly positive on the interior.
\end{remark}

\begin{figure}[h!]
  \centering
  \includegraphics{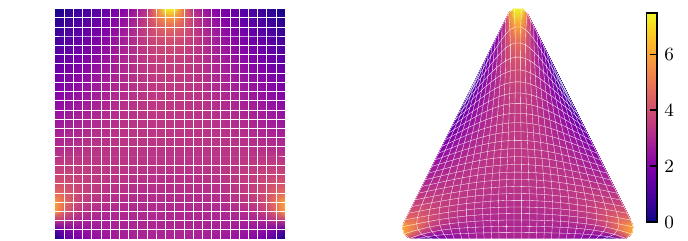}
  \caption{An example square domain (left) mapped to a smooth star domain (right) by solving the
associated Laplacian.  Shading on the domains denotes values of the  Jacobian determinant, which only vanishes on
the boundary, implying that the map is smooth on the interior.}\label{fig:harmonic_triangle}
\end{figure}

With careful construction of the boundary map $\phi$ to ensure unimodality (i.e., boundary
orientation is preserved, $\phi$ is injective, etc.), the map $T^{-1} : \Omega_c \to \Omega_p$ for convex
$\Omega_p$ can be guaranteed to be smooth in two dimensions.  For example, a mapping from a square
to a smoothed triangle can be seen in \cref{fig:harmonic_triangle}.

\subsection{Learned Mappings}\label{subsec:learning}

While harmonic maps can offer smooth, diffeomorphic transformations under suitable conditions, in
practice we consider domains where convexity is either not satisfied or the analytic
boundary map cannot be constructed.  To this end, we extend our framework to more complex domains by
instead learning an approximate harmonic mapping, $T_{\theta} : \Omega_p \to \Omega_c$, parameterized
by some learnable values $\theta$.

In machine learning, Neural Ordinary Differential Equations~\cite{NODE} (Neural ODEs, or
NODEs), are a continuous generalization of feedforward networks wherein the layer update is given by
some continuous function $\vec{f}(\vec{x}, t) : \mathbb{R}^d \times \mathbb{R} \to \mathbb{R}^d$;
evaluation of the neural network is then computed by integrating $f$ in time.  In this context, we
use Neural ODEs as they guarantee a function that is smooth and invertible (under loose
assumptions)---sufficient for use as our domain mapping.

Formally, we have the ODE given by
\begin{align}
  \frac{d\vec{z}}{dt} (t) &= \vec{f}_{\theta}(\vec{z}(t), t), \label{eqn:ode_deriv} \\
  \vec{z}(0) &= \vec{x}, \label{eqn:ode_ic}
\end{align}
where $\vec{f}_{\theta}(\x, t)$ is a learned vector field with parameters $\theta$.  Rewriting this
in integral form, we obtain an expression for the mapping as
\begin{equation}
  T_{\theta}(\x) = \vec{z}(1) = \int_0^1 \vec{f}_{\theta}(\vec{z}(t), t) \,dt
\end{equation}
where $T_{\theta}$ is now a learned map parameterized by $\theta$.  Derivatives of the network
evaluation, for both backpropagation and finite-element assembly, can be efficiently computed by
integration of an adjoint equation, see \cite{Pontryagin1987,NODE} and \cref{sec:hessian_adjoint}.
The inverse of the function is computed by integrating backwards in time.

To train $T_{\theta}$ to behave approximately like a harmonic mapping, we will define a composite
PINNs-style~\cite{RaissiPINNs2017} formulation, following~\cref{thm:diffeo}.  Restating this explicitly,
we would like to satisfy
\begin{align}
  -\Delta T_{\theta} &= \mat{0} \quad \text{in } \Omega_p, \label{eqn:t_int} \\
  T_{\theta} (\partial \Omega_p) &= \Omega_c. \label{eqn:t_bdy}
\end{align}

\begin{figure}[h!]
  \centering
  \includegraphics{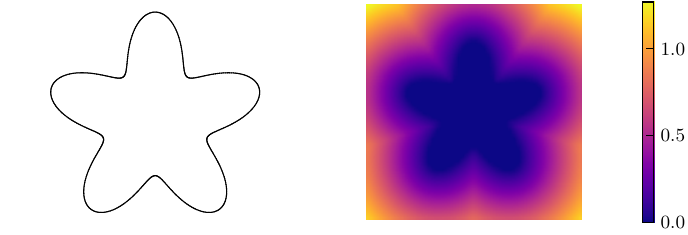}
  \caption{An example star-shaped domain (left) with corresponding unsigned distance function
(right).  Even if the underlying geometry is non-convex it is possible to construct a distance
function to represent it.}\label{fig:df_star}
\end{figure}
First, consider the boundary conditions of the map, which we impose weakly through a loss
function inspired by large deformation diffeomorphic metric mappings
({LDDMMs})~\cite{Hernandez2024}.  Let $I(\cdot; {\Omega})$ be the unsigned distance function defined
as
\begin{align}
  I(\x; \Omega) =
  \begin{cases}
    \|\x - P(\x; \partial \Omega)\|_2 & \x \not \in \Omega \\
    0 & \text{otherwise}
  \end{cases}
\end{align}
where $P(\x; \partial \Omega) = \argmin_{\vec{y} \in \partial \Omega} \|\vec{y} - \x\|^2_2 $ is
the orthogonal projector onto the boundary of $\Omega$.  This is referred to as the \emph{image} of the domain,
and can be interpreted as encoding the interior volume and boundary of the underlying geometry, as
shown in \cref{fig:df_star}.
We emphasize that both of these functions are defined over all $\mathbb{R}^d$, not simply in their
respective domains.

Next, we define the loss by
\begin{equation}
  \ell(\theta) = \|I_p \circ T^{-1} - I_c\|^2_{L^2(\Omega_c)} + \alpha \big\|\Delta T\|^2_{L^2(\Omega_c)}, \label{eqn:lddmm_loss}
\end{equation}
where $I_p(\x) = I(\x,\Omega_p)$ and $I_c(\x) = I(\x,\Omega_c)$, so that the first term weakly enforces the boundary condition of the harmonic map. Intuitively, points
should preserve their distance to the domain boundary after being transformed with $T$.  The second
term in the loss encourages the resulting function to behave approximately like a harmonic: $\Delta
T \in \mathbb{R}^d$ is the componentwise Laplacian of $T$ and $\alpha>0$ is a regularization
parameter penalizing large Laplacian values.  Computation of this componentwise Laplacian is detailed in
\cref{sec:hessian_adjoint}.

The integral in \cref{eqn:lddmm_loss} is approximated by a Monte Carlo scheme of the form
\begin{equation}
  \ell(\theta) \approx \frac{1}{N_s} \sum_{i=1}^{N_s} \left( (I_p (T^{-1}(\x_i)) - I_c(\x_i))^2 +
  \alpha\big\|\Delta T(\x_i)\big\|^2_2\right),
\end{equation}
where the points $\{ \x_i \}$ are drawn uniformly from $\mathcal{U}(-1.5, 1.5)^2$, a region
extending slightly past the computational domain.

One consequence of learning on domains as compared to meshes is that the mapping needs only to be
computed once if the underlying domain remains the same.  For example, if adaptive mesh refinement
or a moving mesh method is used, then $T_{\theta}$ remains a valid mapping as the domain itself does
not change.  This also allows for the training cost to be amortized for, e.g., time-stepped
simulations.

\section{Numerical Results}\label{sec:numerics}

In this section, we compare the convergence of our diffeomorphic multigrid solver to
``off-the-shelf'' Ruge-St\"{u}ben algebraic multigrid (RS-AMG) solver~\cite{Ruge_1987}, as
implemented in {PyAMG}~\cite{pyamg}.  For all of these tests we consider $\mat{A}_p$ to come from
discretizing the physical-grid PDE with piecewise linear elements on triangles for a problem with a
zero right-hand side for $\vec{f}_p$ and a random initial guess for $\vec{u}_p$.  We use stationary
iterations to solve until the Euclidean norm of the residual is reduced below the value $10^{-10}$.

With our diffeomorphic multigrid solver, we now have the additional choice of selecting the size of
the computational mesh as a solver parameter.  This mesh size directly affects the resulting
convergence and work done by the solver: we would like the computational cost of the V-cycle of our
method to be as low as possible, while making the algorithm robust and efficient with the best
possible convergence rate --- these objectives are at odds, as too coarse of a computational mesh will
lead to poor corrections, while a much finer computational mesh will be too costly to run as a
practical method.  In this section, we will fix the computational mesh so that the resulting solver
roughly does the same amount of work as the RS-AMG solver, though we
show the behavior of the method as this mesh size is varied in \cref{subsec:comp_scaling}.

To make the relative work of our solver comparable to that of AMG, we aim to minimize the difference
between the overall \emph{grid complexity} (GC) of the solvers, defined as the sum of unknowns on
all levels divided by the number of unknowns of the finest level,
\begin{equation}
  \text{GC} = \frac{\sum_{\ell=1}^{L} N_{\ell}}{N_{1}},
\end{equation}
where $N_\ell$ is the number of degrees of freedom (DoF) on level $\ell$.  We construct a comparable
diffeomorphic solver by first running the RS-AMG setup on the physical problem to determine its grid complexity,
then experimentally determining a uniform grid size that minimizes the difference of the resulting
grid complexities, favoring higher complexity on the D-MG solver to lower.
We also measure the operator complexity of the solvers as a metric of asymptotic storage cost,
\begin{equation}
  \text{OC} = \frac{\sum_{\ell=1}^{L} \nnz_{\ell}}{\nnz_{1}},
\end{equation}
where $\nnz_\ell$ is the number of nonzeros in the stored operator on level $\ell$.  We note that we do
not use this measure in the setup of our solver, but just to evaluate its efficiency.

\subsection{Varying geometry and diffusivity}

We test our solver against AMG on several domains in two dimensions, including both domains with
known conformal analytic mappings and ones where the mapping must be learned.  These domains and
test problems will be detailed below, as well as their results summarized in
\cref{tab:results_summary1,tab:results_summary2}.  Unless otherwise stated we take $\D = I$.

\textbf{Quarter annulus}
A cross-section of an annulus extending from $\theta = \pi/2$ to $\pi$, with inner radius $r_2=1/2$ and outer
radius $r_1=1$.  This has a known mapping to the square domain given by
\begin{equation}
  T(x, y) =
  \begin{bmatrix}
    2 (\sqrt{x^2 + y^2} - r_2) / (r_1-r_2) - 1\\
    \frac{4}{\pi}\atantwo(y, x) - 3
  \end{bmatrix}.
\end{equation}
Here, we use the two-argument arctangent function, $\atantwo$, that returns $\arctan(y/x)$ when
$x>0$, but properly extends the range of arctangent from $-\pi$ to $\pi$ for negative values of $x$.

\textbf{Restricted channel}
A domain consisting of the unit square, $(-1,1)^2$, with two semicircle ``cut-outs'' of radius
$r=0.1$ located on the top and bottom edges, centered at $(0, 1)$ and $(0, -1)$, respectively.  This can be seen
as a stretched version of the classic restricted channel test problem in fluid dynamics with an
equal aspect ratio.  We learn the mapping, $T$, by a neural ODE.

\textbf{Wavy box}
A deformed, non-convex box with boundary defined by
\begin{align}
  x(\phi) &= \Big( \big( |\cos\phi|^s + |\sin\phi|^s \big)^{-\frac{1}{s}} (W/2) + A\cos(k\phi) \Big) \cos\phi + W/2, \\
  y(\phi) &= \Big( \big( |\cos\phi|^s + |\sin\phi|^s \big)^{-\frac{1}{s}} + A\cos(k\phi) \Big) \sin\phi,
\end{align}
for polar angle $-\pi < \phi \leq \pi$ and parameters $s=500$, $W=3$, $A=0.08$, $k=6$.
While it may be possible to construct a smooth map to the computational domain, we elect to learn it via a neural ODE.

\textbf{Smooth star}
A 5-legged star with smoothed edges, with boundary defined by
\begin{align}
  x(\phi) &= \frac{1}{2}\Big(1 + \cos\Big(\frac{5\phi}{2}\Big)^2\Big)\cos\phi, \\
  y(\phi) &= \frac{1}{2}\Big(1 + \cos\Big(\frac{5\phi}{2}\Big)^2\Big)\sin\phi,
\end{align}
for polar angle $-\pi < \phi \leq \pi$ .  This domain does not have an analytical conformal mapping, nor is it convex; using our methodology, we
learn the mapping $T$ parameterized by a neural ODE.

\textbf{Discontinuous circular coefficients}
A unit square domain with discontinuous diffusivity coefficient given by
\begin{equation}
  \mathcal{D}(x,y) =
  \begin{cases}
    100 & \sqrt{x^2 + y^2} < 3/5 \\
    1 & \text{otherwise}
  \end{cases}.
\end{equation}
The physical domain is meshed to align with the circular discontinuity, while the computational domain is not.
As the physical domain matches our computational domain here, we simply set $T$ to be the identity.

\textbf{Discontinuous checkerboard coefficients}
A unit square domain with rapidly oscillating checkerboard diffusivity coefficients, given by
\begin{equation}
  \mathcal{D}(x,y) =
  \begin{cases}
    100 & (\lfloor 4x \rfloor + \lfloor 4y \rfloor) \equiv 0 \mod 2 \\
    1 & \text{otherwise}
  \end{cases},
\end{equation}
where we use the floor function to define where the coefficient is large.
The physical domain is not meshed to align with the discontinuous coefficients.
Again, because the domains coincide, the identity mapping is used.

\textbf{Discontinuous thin layer coefficients}
A unit square domain with a single vertical ``stripe'' discontinuity, given by
\begin{equation}
  \mathcal{D}(x,y) =
  \begin{cases}
    10^{-10} & 0 \leq x < \varepsilon \\
    2 & \text{otherwise}
  \end{cases},
\end{equation}
for $\varepsilon=0.06$.  The physical domain is not meshed to align with the discontinuous
coefficients.  Again, because the domains coincide, the identity mapping is used.

\begin{table}[!tp]
  \begin{tabular}[t]{cc}
  \toprule
  \multicolumn{2}{c}{\textbf{Quarter annulus}} \\
  \includegraphics{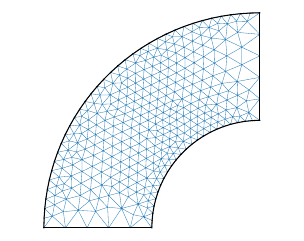} &
  \includegraphics{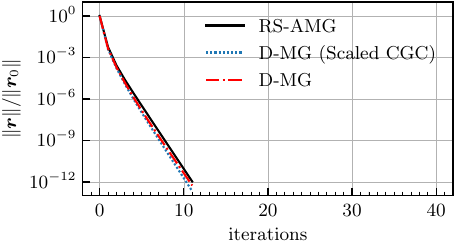} \\
  \midrule
  \multicolumn{2}{c}{\textbf{Restricted channel}} \\
  \includegraphics{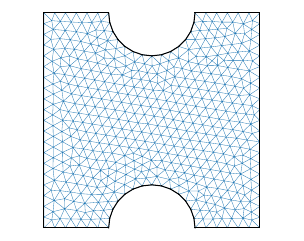} &
  \includegraphics{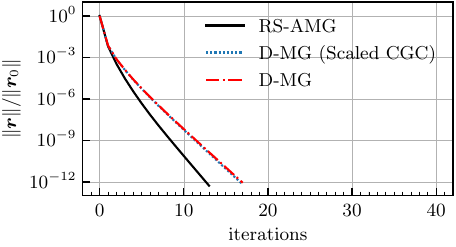} \\
  \midrule
  \multicolumn{2}{c}{\textbf{Wavy box}} \\
  \includegraphics{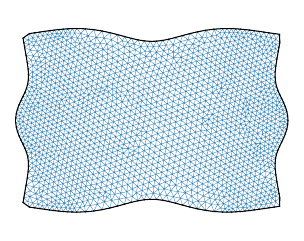} &
  \includegraphics{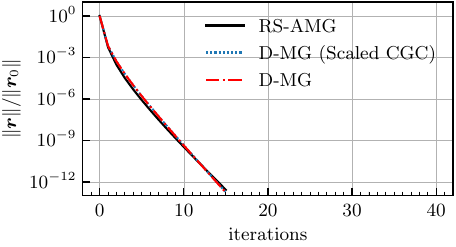} \\
  \midrule
  \multicolumn{2}{c}{\textbf{Smooth star}} \\
  \includegraphics{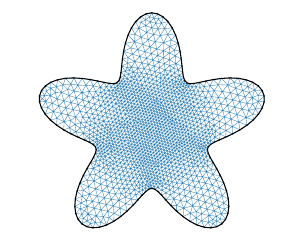} &
  \includegraphics{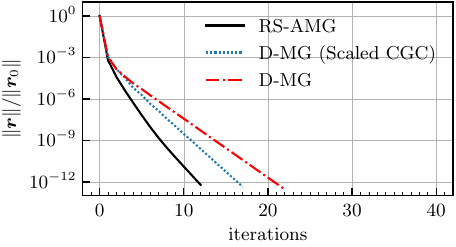} \\
  \bottomrule
  \end{tabular}
  \caption{Physical mesh (left) and residual plot (right).
{RS-AMG} refers to a standard Ruge-St\"{u}ben AMG solver, while {D-MG} and {D-MG}
(scaled CGC) refer to our diffeomorphic multigrid solver, without and with coarse-grid scaling,
respectively.}\label{tab:results_summary1}
\end{table}

\begin{table}[!tp]
  \begin{tabular}[t]{cc}
  \toprule
  \multicolumn{2}{c}{\textbf{Discontinuous circular coefficients}} \\
  \includegraphics{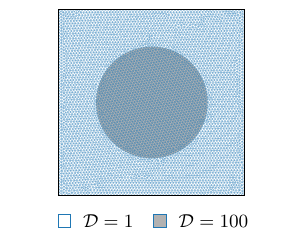} &
  \includegraphics{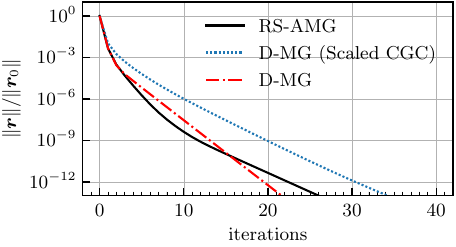} \\
  \midrule
  \multicolumn{2}{c}{\textbf{Discontinuous checkerboard coefficients}} \\
  \includegraphics{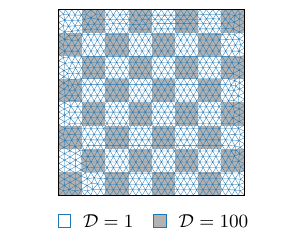} &
  \includegraphics{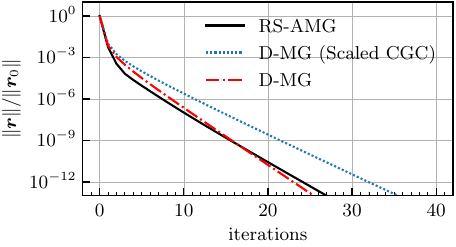} \\
  \midrule
  \multicolumn{2}{c}{\textbf{Discontinuous thin layer coefficients}} \\
  \includegraphics{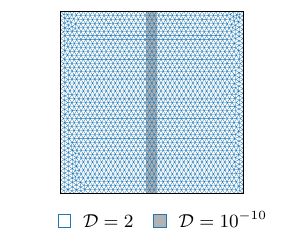} &
  \includegraphics{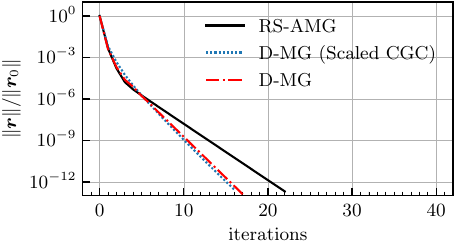} \\
  \bottomrule
  \end{tabular}
  \caption{Physical mesh with diffusivity in gray/white. (left) and residual plot (right).
{RS-AMG} refers to a standard Ruge-St\"{u}ben AMG solver, while {D-MG} and {D-MG}
(scaled CGC) refer to our diffeomorphic multigrid solver, without and with coarse-grid scaling,
respectively.}\label{tab:results_summary2}
\end{table}

\begin{table}[!htp]
  \centering
\begin{tabular}{
  l                     %
  S[table-format = 3.0] %
  @{\hskip 20pt}
  S[table-format = 3.0] %
  S[table-format = 1.2] %
  S[table-format = 1.2] %
  S[table-format = 1.2] %
  @{\hskip 20pt}
  S[table-format = 1.2] %
  S[table-format = 1.2] %
}
  \toprule
   & &
\multicolumn{4}{c}{D-MG} &
\multicolumn{2}{c}{RS-AMG} \\
\cmidrule(r{20pt}){3-6}
\cmidrule(r{4pt}){7-8}
& $N_p$ & $N_c$ & $\gamma$ & GC & OC& GC & OC \\
\cmidrule(r{20pt}){1-2}
\cmidrule(r{20pt}){3-6}
\cmidrule(r{4pt}){7-8}
Quarter annulus                 & 281  & 144  & 0.93 & 1.49 & 1.56 & 1.41 & 1.62 \\
Restricted channel              & 430  & 196  & 0.99 & 1.45 & 1.52 & 1.41 & 1.61 \\
Wavy box                        & 1708 & 576  & 1.00 & 1.38 & 1.46 & 1.41 & 1.69 \\
Smooth star                     & 1381 & 529  & 1.10 & 1.44 & 1.53 & 1.39 & 1.60 \\
Disc. circular coefficients     & 5031 & 1764 & 2.03 & 1.43 & 1.53 & 1.44 & 1.81 \\
Disc. checkerboard coefficients & 938  & 361  & 1.36 & 1.44 & 1.52 & 1.42 & 1.56 \\
Disc. thin layer coefficients   & 1777 & 729  & 1.02 & 1.49 & 1.59 & 1.47 & 1.89 \\
\bottomrule
\end{tabular}
  \caption{Physical degrees of freedom ($N_p$),
  number of degrees of freedom on the finest computational grid ($N_c$),
  coarse-grid scaling parameter ($\gamma$),
  grid complexity (GC), and operator complexity (OC)
  for the Diffeo-MG (D-MG) and Ruge-St\"uben AMG (RS-AMG) solvers.}\label{tab:solver_complexity}
\end{table}

In \cref{tab:results_summary1,tab:results_summary2}, we display residual convergence results for
standard Ruge-St\"{u}ben AMG (denoted {RS}-{AMG}) and compare to our diffeomorphic solver framework
(denoted D-MG).  We also include the residual for a diffeomorphic solver both with and without a
scaled coarse-grid-correction (see \cref{subsec:cgc_scaling}).  For domains for which analytic
conformal mappings exist, the performance of D-MG solver is able to match or out-perform that of the
RS-AMG solver; we see this for the quarter annulus and discontinuous coefficient problems.  The gain
in convergence factor is not particularly significant except in the discontinuous circular
coefficients and discontinuous thin layer problems (5 or so iterations); this is likely due to the
geometric multigrid solver used for the computational solve being robust when applied to anisotropic
problems.  We note that using the scaled coarse-grid correction process generally leads to poorer convergence of the D-MG solver for these variable-coefficient problems.  This can be explained by the
$L^2$ interpolation from the computational to physical grids being only geometric in nature: its setup uses the structure of the basis
functions only; it does not, for example, have any dependency on the operator or diffusivity
coefficient itself.  We leave possible modifications of the interpolation operator to account for such variation in coefficients for future work.

The solver results are somewhat more varied for domains in which a learned mapping is applied.  For the wavy box example, we achieve essentially identical behaviour to RS-AMG, but see slower convergence for D-MG applied to the smooth star and restricted channel domains.  For the restricted channel domain, we lose about 5 iterations to convergence both with and without scaling the coarse-grid correction.  Without scaling the coarse-grid correction, D-MG requires about 10 more iterations for convergence on the smooth star domain, but this is improved to only losing 5 iterations with such scaling.  We propose two possible reasons for this degradation of convergence.  First, these are both cases where we learn the mapping in an approximate sense and, for example, do not guarantee that we preserve key properties such as the area of the two domains, while the deformation for the wavy box case is much more simple (and, perhaps, easier to learn accurately).  To address this, we could adapt the loss function for the neural ODE to emphasize appropriate geometric terms.  Secondly, these are both cases where the most deformation occurs in the mappings to the computational domain and, thus, the Jacobian of these maps is most ill-conditioned.  It is possible we could improve performance here by mapping to structured grids on
domains other than the unit square, or by using domain-decomposition approaches to map sub-meshes to
structured grids, but we do not consider these here.

We include grid and operator complexities for each solver in \cref{tab:solver_complexity}, along
with grid sizes for the physical domains and those of the computational meshes.  As noted above, we
fix the size of the computational mesh for each problem to attain a similar grid complexity for our
solver as is observed for the corresponding AMG solver.  The resulting operator complexities of our
diffeomorphic solvers, however, are always less than those of AMG.  This is due to the bounded
complexity in the Galerkin coarse-grid operators in BoxMG, having at most 9-point stencils in 2D.
We remark that the operator complexities could be reduced further (in both cases) if a symmetric
storage scheme was used, reducing the overhead to only five stored values per degree-of-freedom.

\subsection{Varying size of computational grids}\label{subsec:comp_scaling}

\begin{figure}[!htp] \centering \includegraphics{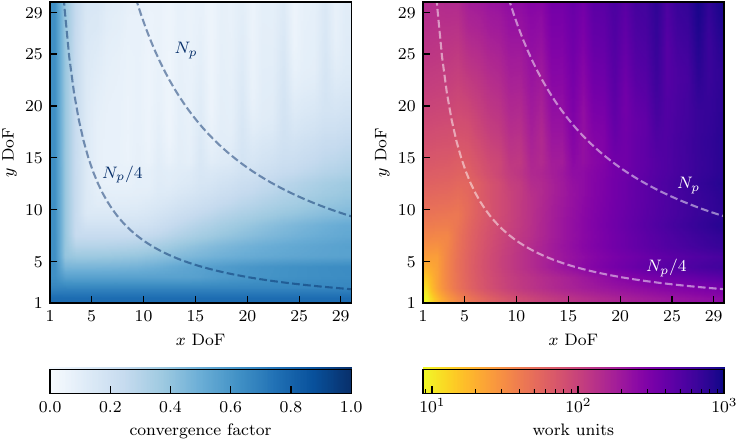}
  \caption{Convergence factor (left) and relative work-per-digit-accuracy (right) of the
D-MG solver on the quarter-annulus domain as the resolution of the computational grid is varied.
$x$-axis indicates width of computational mesh (in elements), while $y$-axis indicates height.  Dashed lines denote
computational meshes with equivalent number of degrees of freedom as the physical mesh ($N_p$) and
coarsening each spatial dimension by two ($N_p/4$).}\label{fig:computational_scaling}
\end{figure}

To study the effect of the resolution of the computational mesh, we selected the quarter annulus
domain and varied the dimensions of the computational mesh, recording the resulting convergence
factor and how it changes as a function of mesh size; the results of this are displayed in
\cref{fig:computational_scaling}, along with the relative work-per-digit-accuracy, given by
\begin{equation}
  \text{wpd}(N_x, N_y) \propto -\frac{N_xN_y}{\log_{10}\rho(N_x, N_y)},
\end{equation}
where $\rho(N_x,N_y)$ is the convergence factor of the D-MG solver with a $N_x \times N_y$
computational mesh.  With this definition, we look for small values of $\text{wpd}(N_x,N_y)$ to reflect efficient solvers.
For this particular domain, while an equal aspect ratio usually suffices to
give a convergent D-MG solver, having higher resolution in the $y$-axis will shift the resulting
solver to have more rapid convergence; this behavior is due to our implementation of the map.  The
elongated portion of the quarter annulus is ``straightened out'' and mapped onto the $y$-axis of the
computational domain.  Hence, higher resolution in the $y$-axis allows the computational mesh to more
accurately represent higher-frequency error modes.  In general, one would expect higher resolution
in both axes to result in better convergence factors, with the trade-off of having more work on the
computational mesh.  There exists a trade-off in convergence factor with work performed, as shown in the plot on the right: for the
quarter annulus domain, best performance is achieved when approximately coarsening by 2 in each spatial direction, as
evidenced by the work-per-digit-accuracy plot.

\subsection{Varying problem size}\label{subsec:phys_scaling}

To demonstrate convergence scalability, we select the quarter annulus domain and vary
the resolution of the physical mesh from $329$ DoF to $\num{1066204}$ DoF.
The computational mesh sizes were selected by
roughly coarsening by a factor of $\alpha=0.28$ with an aspect ratio of $N_y=1.5 N_x$, as observed
in the results in \cref{subsec:comp_scaling}.  This is summarized by the equalities
\begin{equation}
  N_x = \Bigg\lceil \sqrt{\alpha\frac{N_p}{1.5}} \Bigg\rceil, \quad
  N_y = \big\lceil 1.5 N_x \big\rceil.
\end{equation}
We solve this problem using our code, as described above and denoted D-MG below, and using the RS-AMG solver from hypre~\cite{hypre,boomeramg} with default parameters.

\begin{figure}[!htp] \centering \includegraphics{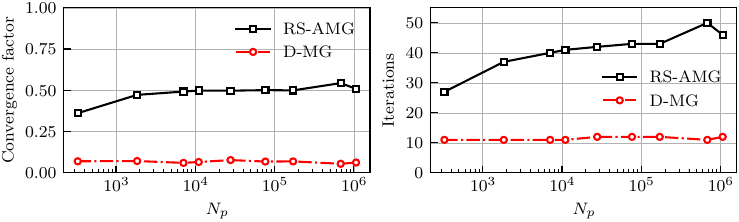}
  \caption{Convergence data for the D-MG and RS-AMG solvers as the problem size increases, indicating that both convergence factor and iterations to convergence remain stable for D-MG, but not RS-AMG as the problem size increases.}\label{fig:convergence}
\end{figure}

From~\cref{fig:convergence}, we see that the D-MG solver maintains consistent convergence factors and iteration counts as the problem size is refined.  The same is not true for RS-AMG, where we see a strong increase in the number of iterations, from 27 for the smallest problem up to 46 for the largest, with corresponding convergence factors increasing from 0.36 to 0.51.  While RS-AMG often achieves much better convergence factors for structured-grid discretizations of the Poisson problem on regular domains, this example shows how that behaviour is somewhat anomalous, at least without substantial parameter tuning.

\begin{figure}[!htp] \centering \includegraphics{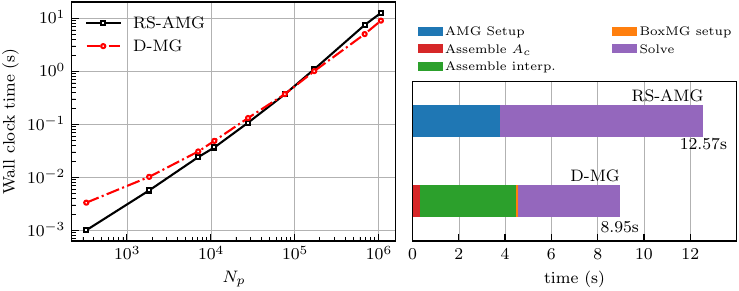}
  \caption{Timing data for the D-MG and RS-AMG solvers.  At left, total time to solution with varying size of the physical mesh.  At right, a breakdown of the timing for the largest problem size. }\label{fig:timings}
\end{figure}

\Cref{fig:timings} presents both recording timings for the scaling study and a breakdown of timings for the largest test case, both recorded running in serial on an Intel Macbook Pro with 2.4 GHz i9 processor and 32GB of DDR4 RAM.  Looking at the solver performance at left, the RS-AMG method tends to lead to slightly faster wall-clock
times for smaller problems, while our D-MG solver has faster runtimes for larger problems
($>10^5$ DoF), with a speedup of about 3 seconds (about 25\%) on the largest problem.  Breaking this down, this can be traced to the slightly longer setup time for the D-MG algorithm, which requires computing the geometric $L^2$ interpolation as
described in \cref{subsec:interpolation}.  While our setup is slightly more involved, having to
compute both the interpolation operator and rediscretize the physical PDE on the computational mesh, this is
counteracted by a faster overall solve time, as seen in the right panel.

It is clear that assembling the interpolation operator is the main computational bottleneck in
the setup of our method, even though we do not invert the mass matrix directly.  While our implementation reasonably uses sparse matrix operations to
compute the interpolation operator, there are still optimizations that could be done to increase the
efficiency of the operator assembly routine.  For example, using some spatial tree data
structure~\cite{BentleyKDTree,FinkelQuadTree} to more efficiently determine overlap of the physical
and computational mesh elements, using a matrix-free approach to applying the operator, or even
restructuring the routine to take advantage of parallelism are potential future directions to speed
up the operator assembly.  However, even with this slightly higher setup cost than
traditional RS-AMG, the resulting speedup in solve times is enough to counteract this on larger problems.  Notably, for the largest problems, we require about 1/4 the number of iterations as RS-AMG to achieve a relatively strict solver tolerance, although our cost per iteration is higher than that of RS-AMG.  This is due to the need to solve with the physical grid mass matrix twice per V-cycle (once for restriction from the physical grid to the finest computational grid, and once for interpolation back to the physical grid), which is done using just two iterations of unpreconditioned CG, which still incurs a similar cost as relaxation on the physical grid.  Further approximating the mass solve, such as via mass lumping, may lead to improvements here.

\section{Conclusions}\label{sec:conclusions}
We present a new framework for enabling the use of geometric multigrid on more complex
geometries via the use of (learned) diffeomorphic mappings to transfer computational work to
structured domains in an auxiliary space multigrid framework.  This not only preserves the efficiency and scalability of robust geometric multigrid
methods, but also allows their use on problems that are more traditionally solved by algebraic
methods.  Our method shows comparable performance to off-the-shelf AMG in terms of number of iterations
to problem convergence.  The diffeomorphic mapping approach allows us to reframe
more difficult problem domains in a multilevel geometric setting, with the inherent benefits therein.

Possible future directions include validation of the use of learned mappings in this solver
framework.  It may very well be the case that there are certain \emph{reference shapes} that are
more well suited to such computations: for example, mapping a circle domain to a square proves
challenging due to the resulting singularities at the corners, but mapping to a different structured
computational grid may be more feasible.  Another interesting direction would be to consider applying this in a
domain decomposition setting: geometric multigrid is already used in practice for semi-structured
meshes, and allowing more flexibility in the types of (sub)-domains would allow such solvers to be
much more competitive with existing algebraic methods in terms of the geometries to which they can be
applied.  Finally, exploring the implementation of this method on accelerators, such as {GPU}s, could
result in significant speedup, both relative to the CPU implementation and algebraic methods
on the {GPU}.

\section*{Acknowledgements}
This material is based in part upon work supported by the Department of Energy, National Nuclear
Security Administration, under Award Number {{DE-NA0003963}}.  The work of SM was partially
supported by an NSERC Discovery Grant.

This research used resources provided by the Darwin testbed at Los Alamos National Laboratory (LANL)
which is funded by the Computational Systems and Software Environments subprogram of LANL's Advanced
Simulation and Computing program (NNSA/DOE).

\appendix
\section{Adjoint computation of ODE Hessian}\label{sec:hessian_adjoint}
We assume we have some ODE solution $\vec{z} : [0,1] \to \mathbb{R}^d$, that is at least $C^2$
with respect to both time and the initial condition ($\x$), and satisfies the ODE: i.e., \cref{eqn:ode_deriv,eqn:ode_ic}.  The Jacobian and Hessian of $\vec{z}$ with respect
to the initial condition, $\x$, can be computed by integrating a system of adjoint equations.

We will first explicitly define the Jacobian matrix and Hessian tensor entrywise, as
\begin{align}
  \big[\mat{J}_z\big]_{ij} &= \bigg[ \frac{\partial \vec{z}}{\partial \x} \bigg]_{ij}
                           = \frac{\partial z_i}{\partial x_j}, \\
  \big[\tens{H}_z\big]_{ijk} &= \bigg[ \frac{\partial^2 \vec{z}}{\partial \x^2} \bigg]_{ijk}
                             = \frac{\partial^2 z_i}{\partial x_j \partial x_k},
\end{align}
where $\mat{J}_z \in \mathbb{R}^{d \times d}$ and $\tens{H}_z \in \mathbb{R}^{d \times d \times d}$
are assumed to be evaluated at some known time, $t$.  Additionally, we will define the adjoint variables
$\mat{W}=\mat{J}_z$ and $\tens{Q}=\tens{H}_z$.  We will show that these can be computed by
integrating an adjoint ODE system.

First, consider the evolution of the Jacobian, $\mat{W}$.  Taking the derivative with respect to
time, we obtain
\begin{align}
  \frac{\partial \mat{W}}{\partial t} &= \frac{\partial}{\partial t} \frac{\partial \vec{z}}{\partial \x} \\
  \intertext{Because of our assumption on smoothness, we can interchange the order of the partial derivatives, giving}
  \frac{\partial \mat{W}}{\partial t} &= \frac{\partial}{\partial \x} \frac{\partial \vec{z}}{\partial t} \\
                                      &= \frac{\partial}{\partial \x} \vec{f}_{\theta}(\vec{z}(t), t) \\
                                      &= \frac{\partial \vec{f}_{\theta}}{\partial \vec{z}} \frac{\partial \vec{z}}{\partial \vec{x}}= \frac{\partial \vec{f}_{\theta}}{\partial \vec{z}} \mat{W}.
\end{align}
The initial condition is given simply by $\mat{W}(0) = \mat{I}$, since $(\partial \vec{z}/\partial
\x)(0)=\partial \x/\partial \x$.

A similar derivation can be done to find the time evolution of the Hessian, $\tens{Q}$.  Again,
taking the derivative with respect to time, we obtain
\begin{align}
  \frac{\partial \tens{Q}}{\partial t} &= \frac{\partial}{\partial t} \frac{\partial^2 \vec{z}}{\partial \x^2} \\
                                       &= \frac{\partial^2}{\partial \x^2} \frac{\partial \vec{z}}{\partial t} \\
                                       &= \frac{\partial^2}{\partial \x^2} \vec{f}_{\theta}(\vec{z}(t), t). \label{eqn:w_dt}
\end{align}
Applying the second-order multivariate chain rule as given by Fa\`{a} di Bruno's
formula~\cite{Fraenkel1978}, we obtain the expression
\begin{equation}
  \bigg[ \frac{\partial \tens{Q}}{\partial t} \bigg]_{ijk} =
  \sum_{l=1}^d \bigg( \frac{\partial f_i}{\partial z_l} q_{ljk} \bigg) +
  \sum_{l=1}^d\sum_{m=1}^d \bigg( \frac{\partial^2 f_i}{\partial z_l \partial z_m} w_{lj} w_{mk} \bigg). \label{eqn:q_dt}
\end{equation}
The initial condition is given by the zero tensor, $\tens{Q}(0)_{ijk} = 0$.

Computing $\mat{W}$ and $\tens{Q}$ for some time $t$, can then be done by integrating the coupled
ODE system $(\mat{W}, \tens{Q})$ using \cref{eqn:w_dt,eqn:q_dt}.  The componentwise Laplacian of
$\vec{z}$ can then be computed by taking traces of horizontal slices of $\tens{Q}$,
\begin{equation}
  \big[\Delta \vec{z}(t)\big]_i = \tr \big(\tens{Q}_{i,:,:}(t)\big).
\end{equation}
Note that element $q_{ijk}$ of the stored Hessian depends only on elements $q_{ljk}$ for $l=1\ldots d$,
meaning that to compute the Laplacian we, in fact, only need to store and evolve the diagonal elements of the
last two indices ($q_{ijj}$ for $i,j=1\ldots d$).

\bibliographystyle{etna}
\bibliography{paper_diffeomg}

\begin{thebibliography}{10}

\bibitem{Alessandrini2001}
{\sc G.~Alessandrini and V.~Nesi}, {\em Univalent $\sigma$-harmonic mappings},
  Archive for Rational Mechanics and Analysis, 158 (2001), pp.~155--171.
\newblock \\ \url{http://link.springer.com/10.1007/PL00004242}.

\bibitem{Bell2012}
{\sc N.~Bell, S.~Dalton, and L.~N. Olson}, {\em Exposing fine-grained
  parallelism in algebraic multigrid methods}, SIAM Journal on Scientific
  Computing, 34 (2012).

\bibitem{pyamg}
{\sc N.~Bell, L.~N. Olson, J.~Schroder, and B.~Southworth}, {\em {PyAMG}:
  Algebraic multigrid solvers in python}, Journal of Open Source Software, 8
  (2023), p.~5495.
\newblock \\ \url{https://doi.org/10.21105/joss.05495}.

\bibitem{BentleyKDTree}
{\sc J.~L. Bentley}, {\em Multidimensional binary search trees used for
  associative searching}, Commun. ACM, 18 (1975), p.~509–517.
\newblock \\ \url{https://doi.org/10.1145/361002.361007}.

\bibitem{BrambleNonnested1991}
{\sc J.~H. Bramble, J.~E. Pasciak, and J.~Xu}, {\em The analysis of multigrid
  algorithms with nonnested spaces or noninherited quadratic forms},
  Mathematics of Computation, 56 (1991), p.~1–34.
\newblock \\ \url{http://dx.doi.org/10.1090/S0025-5718-1991-1052086-4}.

\bibitem{MBrezina_etal_2006a}
{\sc M.~Brezina, R.~Falgout, S.~MacLachlan, T.~Manteuffel, S.~McCormick, and
  J.~Ruge}, {\em Adaptive algebraic multigrid}, SIAM J. Sci. Comput., 27
  (2006), pp.~1261--1286.

\bibitem{BurmanCutFEM2015}
{\sc E.~Burman, S.~Claus, P.~Hansbo, M.~G. Larson, and A.~Massing}, {\em
  Cutfem: Discretizing geometry and partial differential equations},
  International Journal for Numerical Methods in Engineering, 104 (2015),
  pp.~472--501.
\newblock \\ \url{https://onlinelibrary.wiley.com/doi/abs/10.1002/nme.4823}.

\bibitem{ChenAuxiliary2015}
{\sc L.~Chen, J.~Wang, Y.~Wang, and X.~Ye}, {\em An auxiliary space multigrid
  preconditioner for the weak galerkin method}, Computers \& Mathematics with
  Applications, 70 (2015), pp.~330--344.
\newblock \\
  \url{https://www.sciencedirect.com/science/article/pii/S0898122115001972}.

\bibitem{NODE}
{\sc R.~T.~Q. Chen, Y.~Rubanova, J.~Bettencourt, and D.~Duvenaud}, {\em Neural
  ordinary differential equations}, NIPs, 109 (2018), pp.~31--60.
\newblock \\ \url{https://arxiv.org/abs/1806.07366v5}.

\bibitem{ChesshireOverlapping1990}
{\sc G.~Chesshire and W.~Henshaw}, {\em Composite overlapping meshes for the
  solution of partial differential equations}, Journal of Computational
  Physics, 90 (1990), pp.~1--64.
\newblock \\
  \url{https://www.sciencedirect.com/science/article/pii/0021999190901968}.

\bibitem{DendyBoxMG1982}
{\sc J.~Dendy}, {\em Black box multigrid}, Journal of Computational Physics, 48
  (1982), pp.~366--386.
\newblock \\
  \url{https://linkinghub.elsevier.com/retrieve/pii/0021999182900572}.

\bibitem{DendyBoxMG2010}
{\sc J.~E. Dendy and J.~D. Moulton}, {\em Black box multigrid with coarsening
  by a factor of three}, Numerical Linear Algebra with Applications, 17 (2010),
  pp.~577--598.
\newblock \\ \url{https://onlinelibrary.wiley.com/doi/10.1002/nla.705}.

\bibitem{FinkelQuadTree}
{\sc R.~A. Finkel and J.~L. Bentley}, {\em Quad trees a data structure for
  retrieval on composite keys}, Acta Informatica, 4 (1974), pp.~1--9.
\newblock \\ \url{https://doi.org/10.1007/BF00288933}.

\bibitem{Fraenkel1978}
{\sc L.~E. Fraenkel}, {\em Formulae for high derivatives of composite
  functions}, Mathematical Proceedings of the Cambridge Philosophical Society,
  83 (1978), pp.~159--165.
\newblock \\
  \url{https://www.cambridge.org/core/product/identifier/S0305004100054402/type/journal_article}.

\bibitem{GelfandCalculus2000}
{\sc I.~M. Gelfand and S.~V. Fomin}, {\em Calculus of Variations}, Dover Books
  on Mathematics, Dover Publications, Mineola, NY, Oct. 2000.

\bibitem{HenshawOverlapping2005}
{\sc W.~D. Henshaw}, {\em On multigrid for overlapping grids}, SIAM Journal on
  Scientific Computing, 26 (2005), p.~1547–1572.
\newblock \\ \url{http://dx.doi.org/10.1137/040603735}.

\bibitem{HenshawOverlapping2003}
{\sc W.~D. Henshaw and D.~W. Schwendeman}, {\em An adaptive numerical scheme
  for high-speed reactive flow on overlapping grids}, Journal of Computational
  Physics, 191 (2003), pp.~420--447.
\newblock \\
  \url{https://www.sciencedirect.com/science/article/pii/S0021999103003231}.

\bibitem{boomeramg}
{\sc V.~E. Henson and U.~M. Yang}, {\em Boomeramg: A parallel algebraic
  multigrid solver and preconditioner}, Applied Numerical Mathematics, 41
  (2002), p.~155–177.
\newblock \\ \url{http://dx.doi.org/10.1016/S0168-9274(01)00115-5}.

\bibitem{Hernandez2024}
{\sc M.~Hernandez and U.~R. Julvez}, {\em Insights into traditional large
  deformation diffeomorphic metric mapping and unsupervised deep-learning for
  diffeomorphic registration and their evaluation}, Computers in Biology and
  Medicine, 178 (2024).

\bibitem{hypre}
{\em {\sl hypre}: High performance preconditioners}.
\newblock \url{https://llnl.gov/casc/hypre},
  \url{https://github.com/hypre-space/hypre}.

\bibitem{XiangminNonmatching2004}
{\sc X.~Jiao and M.~T. Heath}, {\em Common-refinement-based data transfer
  between non-matching meshes in multiphysics simulations}, International
  Journal for Numerical Methods in Engineering, 61 (2004), pp.~2402--2427.
\newblock \\ \url{https://onlinelibrary.wiley.com/doi/abs/10.1002/nme.1147}.

\bibitem{JuntunenNitsche2009}
{\sc M.~Juntunen and R.~Stenberg}, {\em Nitsche’s method for general boundary
  conditions}, Mathematics of Computation, 78 (2009), pp.~1353--1374.
\newblock \\
  \url{http://www.ams.org/jourcgi/jour-getitem?pii=S0025-5718-08-02183-2}.

\bibitem{LiFNO2023}
{\sc Z.~Li, D.~Z. Huang, B.~Liu, and A.~Anandkumar}, {\em Fourier neural
  operator with learned deformations for pdes on general geometries}, J. Mach.
  Learn. Res., 24 (2023).

\bibitem{LiFNO2021}
{\sc Z.~Li, N.~Kovachki, K.~Azizzadenesheli, B.~Liu, K.~Bhattacharya,
  A.~Stuart, and A.~Anandkumar}, {\em Fourier neural operator for parametric
  partial differential equations}, 2021.

\bibitem{PontInterpolation2016}
{\sc A.~Pont, R.~Codina, and J.~Baiges}, {\em Interpolation with restrictions
  between finite element meshes for flow problems in an ale setting:
  Interpolation with restrictions}, International Journal for Numerical Methods
  in Engineering, 110 (2016), p.~1203–1226.
\newblock \\ \url{http://dx.doi.org/10.1002/nme.5444}.

\bibitem{Pontryagin1987}
{\sc L.~S. Pontryagin}, {\em Mathematical theory of optimal processes:
  Mathematical theory of optimal processes L.s. pontryagin selected works
  volume 4}, Classics of Soviet Mathematics, Harwood Academic, Amsterdam,
  Netherlands, Mar. 1987.

\bibitem{RaissiPINNs2017}
{\sc M.~Raissi, P.~Perdikaris, and G.~E. Karniadakis}, {\em Physics informed
  deep learning (part i): Data-driven solutions of nonlinear partial
  differential equations},  (2017).
\newblock \\ \url{http://arxiv.org/abs/1711.10561}.

\bibitem{ReisnerBoxMG2018}
{\sc A.~Reisner, L.~N. Olson, and J.~D. Moulton}, {\em Scaling structured
  multigrid to 500k+ cores through coarse-grid redistribution},  (2018).

\bibitem{Ruge_1987}
{\sc J.~W. Ruge and K.~Stüben}, {\em 4. Algebraic Multigrid}, Society for
  Industrial and Applied Mathematics, Jan. 1987, p.~73–130.

\bibitem{SaadEigs2011}
{\sc Y.~Saad}, {\em Numerical Methods for Large Eigenvalue Problems}, Society
  for Industrial and Applied Mathematics, 2011.

\bibitem{KStuben_2001a}
{\sc K.~St{\"u}ben}, {\em An introduction to algebraic multigrid}, in
  Multigrid, U.~Trottenberg, C.~Oosterlee, and A.~Sch{\"u}ller, eds., Academic
  Press, London, 2001, pp.~413--528.

\bibitem{XuAuxiliary1996}
{\sc J.~Xu}, {\em The auxiliary space method and optimal multigrid
  preconditioning techniques for unstructured grids}, Computing, 56 (1996),
  pp.~215--235.
\newblock \\ \url{https://doi.org/10.1007/BF02238513}.

\end{thebibliography}

\end{document}